\documentclass[11pt]{article}

\usepackage{latexsym}
\usepackage{amsmath}

 \usepackage[textwidth=6.4in, textheight=8.775in]{geometry}

\newtheorem{theorem}{Theorem}
\newtheorem{lemma}{Lemma}
\newtheorem{prop}{Proposition}
\newtheorem{claim}{Claim}
\newtheorem{observation}{Observation}

\newtheorem{definition}{Definition}

\newcommand{\QED}{$\Box$}

\newcommand{\smallqed}{{\tiny ($\Box$)}}
\newcommand{\proof}{\noindent\textbf{Proof. }}

\newcommand{\cL}{\mathcal{L}}

\let\oldenumerate\enumerate
\renewcommand{\enumerate}{
 	\oldenumerate
\setlength{\itemsep}{1pt}
\setlength{\parskip}{0pt}
\setlength{\parsep}{0pt}
}

\newcommand{\residual}{G_D^\iota}
\newcommand{\NEWresidual}{G_{D'}^\iota}
\newcommand{\expression}{\mathcal{F}}

\newcommand{\w}{\textsf{w}}

\begin{document}

\title{On the isolation number of graphs with minimum degree four}
 	 	
\author{$^{1,2}$Wayne Goddard \, and \, $^2$Michael A. Henning \\ \\
$^1$School of Mathematical and Statistical Sciences \\
Clemson University, United States of America \\
\small \tt Email: goddard@clemson.edu  \\
\\
$^2$Department of Mathematics and Applied Mathematics \\
University of Johannesburg, South Africa\\
\small \tt Email: mahenning@uj.ac.za  \\
}
 	
\date{}
\maketitle

\begin{abstract}
An isolating set in a graph $G$ is a set $S$ of vertices such that removing $S$ and its neighborhood leaves no edge. The isolation number $\iota(G)$ of $G$ (also known as the vertex-edge domination number) is the minimum size among all isolating sets of $G$. We provide a technique for proving upper bounds on this parameter for graphs with a given minimum degree. For example, we show that if $G$ has order~$n$ and minimum degree at least~$4$, then $\iota(G) \le 13n/41$,
and if $G$ is also triangle-free, then $\iota(G) \le 3n/10$.
\end{abstract}

 	
 	
\section{Introduction}
 	
As defined by Caro and Hansberg~\cite{CH},
an isolating set in a graph $G$ is a set $S$ of vertices such that removing $S$ and its neighborhood leaves no edge. The isolation number $\iota(G)$ of $G$  is the minimum size among all isolating sets of $G$. This parameter is equivalent to the vertex-edge domination number, introduced earlier by~Peters~\cite{Peters}. A fundamental upper bound on the parameter was provided by Caro and Hansberg and by \.{Z}yli\'nski: 

\begin{prop} \cite{CH,Zylinski} \label{prop}
If $G$ is a connected graph of order $n$ other than $K_2$ and~$C_5$, then $\iota(G) \le n/3$. 
\end{prop}

The extremal trees were characterized by Krishnakumari et al.~\cite{KVK} and Dapena et al.~\cite{DLSV}, and the extremal graphs were characterized in~\cite{first}. Apart from one graph, all the extremal graphs have a vertex of degree two or less, and there are infinitely many extremal graphs with minimum degree $2$. 

A natural question then is what happens for larger minimum degree. 
In general, Caro and Hansberg~\cite{CH}  showed that the isolation number of a graph of order $n$  is at most $( \ln ( \delta + 1) + \frac{1}{2} ) n / (\delta + 1) $ if every vertex
has degree at least $\delta$. Also, they showed there
are arbitrarily large connected graphs with minimum degree $\delta$ and isolation number $ 2 n / (\delta + 4)$.
Another result in this direction was provided by Ziemann and \.{Z}yli\'{n}ski~\cite{ZZ}, who showed that apart from one graph of order $6$, the isolation number of a cubic graph is at most $10n/31$. They also showed that there are arbitrarily large cubic graphs with isolation number $2n/7$. 

In this paper we introduce a tool called the isolation residual graph, modeled after the residual graph used in bounds on the game domination number (see for example~\cite{BKR,KWZ}) and the domination number  (see for example~\cite{Bu21,BuHe21,BuHe25,BuKl16}).
Using this tool we provide several upper bounds for the isolation number of a graph given a lower bound on its minimum degree.
For example, we show that if $G$ is a graph of order~$n$ with minimum degree at least~$4$, then $\iota(G) \le 13n/41$,
and if $G$ is also triangle-free, then $\iota(G) \le 3n/10$. 

We proceed by first introducing the isolation residual graph and then proving the upper bounds.
Thereafter we briefly discuss associated lower bounds for these problems.

\section{Methodology}

We establish upper bounds on the isolation number of a graph with given minimum degree.
For notation and graph theory terminology, we generally follow~\cite{HaHeHe-23}. 
Specifically, let $G$ be a graph with vertex set $V(G)$ and edge set $E(G)$. For a set of vertices $S\subseteq V(G)$, the subgraph induced by $S$ is denoted by $G[S]$. 

Our proof is based on the proof technique employed to establish upper bounds on the game domination number or domination number of a graph. 
That proof uses a \emph{residual graph}, which is defined by coloring and weighting the vertices of the graph as to how they relate
to a given partial dominating set. We extend this:
 	
\begin{definition}
\label{residual-graph}
\normalfont
Given a graph $G$ and a subset $D \subseteq V(G)$, the \emph{isolation residual graph} 
denoted by $\residual$ is obtained as follows. \\ [-26pt]
\begin{enumerate}
\item 
A vertex in $G-N[D]$ with a neighbor in $G-N[D]$  is colored \emph{white}.  
\item 
A vertex in $N[D]$ with a white neighbor is colored \emph{blue}.
\item
All other vertices are colored \emph{red}.
\item 
The edge set of $\residual$ consists of all edges of $G$ that are incident with at least one white vertex.
\end{enumerate}
\end{definition}

\noindent
In other words, the white vertices still need to be isolated, the blue vertices might be useful to do this isolating, and the red vertices are no longer
relevant.
We denote the set of red vertices by~$R$, the set of blue vertices by $B$, and the set of white vertices by $W$.
We observe that these sets form a (weak) partition of $V(G)$.  	
Note that:

\begin{observation}
\label{ob-residual}
The following holds in the isolation residual graph $\residual$:\\ [-24pt]
\begin{enumerate}
\item
Every vertex of $D$ is colored red.
\item 
Every red vertex is isolated, but no white nor blue vertex is isolated.
\item
The degree of a white vertex in $\residual$ equals its degree in $G$.
\item
The set $D$ is an isolating set in $G$ if and only if there is no white vertex.
\end{enumerate}
\end{observation}
 	
We now introduce a weight function.
We define $B_4$ as the set of blue vertices of degree at least~$4$ in $\residual$, 
and we define $B_i$ as the set of blue vertices of degree exactly~$i$ in $\residual$ for $i \in [3]$. 
The weight $\w(v)$ of a vertex $v$ is defined as:
\[
 \w(v) = \left\{ \begin{array}{l}
  \mbox{$\omega$ if $v\in W$,} \\[1mm]
  \mbox{$\beta_i$ if $v \in B_i$ for $i \in [4]$,  and} \\[1mm]
  \mbox{$0$ if $v \in R$,}
  \end{array}
 \right.
\]
where $\omega$ and each $ \beta_i$ are positive constants to be determined later.
We define the weight of the isolation residual graph $\residual$ as
\[
\w( \residual ) = \sum_{v \in V} \w(v) .
\]
By Observation~\ref{ob-residual}, $\w(  \residual ) = 0$ if and only if $D$ is an isolating set in $G$. 

The proof strategy is to start with $D=\emptyset$ and repeatedly enlarge the set $D$.
We need the following immediate observation.

\begin{observation}
Suppose $D \subseteq D' \subseteq V(G)$.
Every red vertex in $\residual$ is red in
$\NEWresidual$. Every blue vertex in $\residual$ is red or blue 
in $\NEWresidual$.
\end{observation}

For set $A \subseteq V(G) \setminus D$ and $D' = D \cup A$,
define
\[
   \xi (A) = \w( \residual ) - \w( \NEWresidual )
\]
as the decrease of the weight when extending $D$ to $D'$. 
We say that such a set $A$ is a \emph{$D$-desirable set} if
\[
   \xi (A) \ge |A|.
\]
We say a graph $G$ is \emph{alluring} if whenever
$D \subseteq V(G)$ is not an isolating set of $G$,
the graph $G$ contains a $D$-desirable set $A$.

\begin{theorem} \label{t:prototype}
If there is a setting of the weights such that graph $G$ is alluring, then $\iota(G) \le \omega n$ where
$n$ is the order of $G$.
\end{theorem}
\proof
Let $D_0 = \emptyset$. In the isolation residual graph $G_{D_0}^{\iota}$, every vertex is white 
and hence $\w(G_{D_0}^{\iota}) = \omega \, n$. Since $G$ is alluring, there exists a $D_0$-desirable set $A_1$. 
Let $D_1 =A_1$. If $D_1$ is not an isolating set of $G$, then 
there exists a $D_1$-desirable set $A_2$. Let 
$D_2 = A_1 \cup A_2$. If $D_2$ is not an isolating set of $G$, then
we continue the process. 
At the end we obtain an isolating set $D = A_1 \cup \cdots \cup A_k$ of $G$ such that
\[
0 = \w(G_{D_k}^{\iota}) 
\le 
\w(G_{D_{k-1}}^{\iota}) -  |A_k| 
\le \cdots \le 
 \w(G_{D_0}^{\iota} ) -  \textstyle{\sum_{i=1}^k |A_i| } 
= 
\omega n -  |D|.
\]
Thus  $\iota (G) \le |D| \le \omega n$.~\QED

Thus to provide an upper bound for a family of graphs, it suffices to find a setting of the weights such that every graph in 
that family is alluring.


\section{Alluring Graphs}

 For convenience we define 
 \[
 \varepsilon_i = \beta_i - \beta_{i-1}  \quad \mbox{for $2\le i \le 4$.}
 \]
Consider the following list $\cL$ of constraints:

\begin{quote}
\begin{tabular}{ll}
$\omega + 5(\omega-\beta_4) \ge 1$ & (Claim 1) \\
$\omega + 4(\omega-\beta_3) \ge 1$  \\
$\beta_4 + 5(\omega-\beta_3)  \ge 1$ \\[2pt]
$\omega + 3(\omega-\beta_2)  + 4( \delta -3 ) \varepsilon_4 \ge 1$ & (Claim 2)  \\
$\beta_4 + 4(\omega-\beta_2)  + 4( \delta -3 ) \varepsilon_4 \ge 1$ \\[2pt]
$\omega + (\delta-2) \varepsilon_3 \ge \frac{1}{3} $ & (Claim 3) \\[2pt]
$4 \omega + \beta_2  \ge 1$ & (Claim 4) \\
$7 \omega  + \beta_2 \ge 2$ \\
$10 \omega  + \beta_2 \ge 3$ \\[2pt]
$2 \omega + 2 (\delta-1) \min ( \beta_1, \beta_2/2 ) \ge 1$ & (Claim 5) \\
$5  \, \omega + 5 ( \delta - 2)  \min ( \beta_1, \beta_2/2, \beta_3/3) \ge 2$ \\[2pt]
$\omega \ge \beta_4 \ge \beta_3 \ge \beta_2 \ge  \beta_1 > 0$ &  (non-increasing $\beta_i$) \\[2pt]
$ \varepsilon_4 \le \varepsilon_3 \le \varepsilon_2 \le \beta_1$ & (non-decreasing $\varepsilon_i$) 
\end{tabular}
\end{quote}

We prove:

\begin{lemma}  \label{l:alluring}
Assume the weights satisfy the system of constraints $\cL$ for some $\delta \ge 3$. Then every 
graph with minimum degree at least $\delta $ is alluring.
\end{lemma}

We prove Lemma~\ref{l:alluring} by a series of claims, where each claim shows that,
 if some condition holds, a $D$-desirable set is guaranteed to exist.

If sets $D$ and $A$ are clear from the context, we set $D' = A \cup D$ and denote by $W'$, $B'$, and $R'$, 
the set of white, blue, and red, vertices, respectively, in $\NEWresidual$. Similarly 
we define $B_4'$ as the set of blue vertices of degree at least~$4$ in $\NEWresidual$, and  
$B_i'$ as the set of blue vertices of degree exactly~$i$ in $\NEWresidual$ for $i \in [3]$. 

Let $\Delta_W$ denote
the maximum number of white neighbors of a white vertex in $\residual$, and 
$\Delta_B$ denote
the maximum degree of a blue vertex.


\begin{claim}
\label{c:big}
We may assume $\Delta_W \le 3$ and $\Delta_B \le 4$.
\end{claim}
\proof 
(i) Assume $\Delta_W \ge 4$. Let $A = \{ v \}$ where $v$ is a white vertex with $\Delta_W$ white neighbors.
In $\NEWresidual$ the vertex~$v$ is recolored red. Every white neighbor $u$ of $v$ is recolored blue or red,
and has at most $\Delta_W-1$  white neighbors in $\NEWresidual$. 
Thus $u \in R' $ or $u \in B_i'$ for some $i \in [ \min( \Delta_W-1, 4 ) ]$.
Hence  
$ \xi (A) \ge  \expression_{1a} :=  \omega + 5  (\omega-\beta_4) $ 
if $\Delta_W \ge 5$, and
$ \xi (A) \ge \expression_{1b} := \omega + 4  ( \omega - \beta_3)$ 
if $\Delta_W = 4$.
The system $\cL$ ensures $\expression_{1a} , \expression_{1b} \ge 1$, and so
the set $A$ is a $D$-desirable set in each case. 

(ii) Assume $\Delta_B \ge 5$.
Let $A = \{w\}$ where $w$ a blue vertex with $\Delta_B$ (white) neighbors. 
We note that $w \in B_4$ and $w \in R'$; so the weight decrease of~$w$ is~$\beta_4$. 
Let $x$ be an arbitrary (white) neighbor of $w$ in~$\residual$. 
Since by part (i) we may assume $\Delta_W \le 3$, 
we have $x \in R'$ or $x \in B_i'$ with $i \in [3]$ in $\NEWresidual$.
Hence $\xi (A) \ge \expression_{1c} :=  \beta_4 + 5  ( \omega - \beta_3)  $.
The system $\cL$ ensures $\expression_{1c} \ge 1$,
and so the set $A$ is a $D$-desirable set.~\smallqed

It follows that we may assume that a vertex in $B_4$ has degree exactly $4$.
Hence: if $y \in B_i$ has an incident edge $e$ in $\residual$ (necessarily joining it to a white vertex) and edge~$e$ is removed,  
then $y$ has exactly $i-1$ white neighbors. Thus the decrease in weight of $y$ is at least~$\varepsilon_i$.
We note also that we require that the $\varepsilon_i$ are non-decreasing. 

\begin{claim}
\label{c:small}
We may assume $\Delta_W \le 2$ and $\Delta_B \le 3$.
\end{claim}
\proof
(i) By Claim~\ref{c:big} we may assume $\Delta_W \le 3$.
Assume $\Delta_W = 3$.
Let $A = \{ v \}$ where $v$ is a white vertex of with $\Delta_W$ white neighbors. 
In $\NEWresidual$ the vertex~$v$ is recolored red. Every white neighbor $u$ of $v$ is recolored blue or red,
and has at most $\Delta_W-1$  white neighbors in $\NEWresidual$. 
Thus $u \in R' $ or $u \in B_i'$ for some $i \in [ 2 ]$.
By the minimum degree, every white vertex has at least $\delta-3$ edges to~$B$.
Thus there are at least $4 ( \delta-3)$ edges joining $B$ to $v$ and its white neighbors. 
The recoloring causes each of these edges to be removed. 
Hence $\xi (A) \ge  \expression_{2a} := \omega + 3  ( \omega - \beta_2) + 4 ( \delta-3) \varepsilon_4$
in this case. The system $\cL$ ensures $\expression_{2a} \ge 1$, and so
the set $A$ is a $D$-desirable set. 

(ii) By Claim~\ref{c:big} we may assume $\Delta_B \le 4$.
Assume $\Delta_B = 4$.
Let $A = \{w\}$ where 
$w$ is a blue vertex of degree $\Delta_B$. 
We note that $w \in B_4$ and $w \in R'$; so the weight decrease of~$w$ is~$\beta_4$. 
Let $x$ be an arbitrary (white) neighbor of $w$ in~$\residual$. 
Since by part (i) we may assume $\Delta_W \le 2$, 
we have $x \in R'$ or $x \in B_i'$ with $i \in [2]$ in $\NEWresidual$.
Further, vertex $x$ has at least $\delta-3$ blue neighbors different from~$w$. 
Thus, there are at least $\Delta_B( \delta-3)$ edges joining $B \setminus \{ w \}$ to the white neighbors of $w$. 
The removal of each such edge decreases the total weight of $B \setminus \{ w \}$ 
by at least $\varepsilon_4$. 
Hence $\xi (A) \ge \expression_{2b} :=  \beta_4 + 4  ( \omega - \beta_2)  + 4 (\delta-3) \varepsilon_4$.
The system $\cL$ ensures $\expression_{2b} \ge 1$,
and so the set $A$ is a $D$-desirable set.~\smallqed


\begin{claim}
\label{c:whiteComponent}
We may assume every component $F$ in $\residual [W]$ is either $K_2$ or $C_5$.
\end{claim}
\proof 
Assume there is a component $F$ in $\residual [W]$ that is neither $K_2$ nor $C_5$.
Let $A$ be a minimum isolating set of $F$. By Proposition~\ref{prop} it holds that $|A| \le |V(F)|/3$.
(Actually, because we may assume $F$ has maximum degree at most $2$, we need  this fact only for paths and cycles,
and that was shown already by Peters~\cite{Peters}.)
Since we may assume $\Delta_W \le 2$ by Claim~\ref{c:small}, 
there are at least $(\delta-2) |V(F)|$ edges joining $F$ to $B$.
The removal of each such edge decreases the weight of $B$ by at least~$\varepsilon_3$, since by Claim~\ref{c:small}
we may assume $B_4$ is empty. 
Every vertex in $F$ becomes red. Hence $\xi (A) \ge  |V(F)| \times \expression_3$ where $\expression_3 :=   \omega + (\delta-2) \varepsilon_3$.
The system $\cL$ ensures $\expression_3 \ge 1/3$, 
and so $\xi(A) \ge |V(F)|/3 \ge |A|$, so that $A$ is a $D$-desirable set.~\smallqed

\begin{claim}
\label{c:blueTwo}
We may assume there is no blue vertex with neighbors in more than one component in~$\residual[W]$.
\end{claim}
\proof 
Assume blue vertex $x$  is adjacent to white vertices from two components 
in $\residual [W]$, say $F_1$ and $F_2$.
 By Claim~\ref{c:whiteComponent}, for $i \in [2]$ we may assume each $F_i$  is $ K_2$ or $C_5$.
 
Let $y_i$ be a (white) vertex in $F_i$ adjacent to $x$ for $i \in [2]$;
and for each component $F_i$ that is a $5$-cycle, let $u_i$ be a vertex of $F_i$ that is not next to $y_i$ on the cycle.
If both $F_1$ and $F_2$ are $K_2$ let $A = \{ x \}$;
if $F_1=K_2$ and $F_2=C_5$ let $A = \{ x, u_2 \}$; and if
both $F_1$ and $F_2$ are $C_5$ let $A = \{ x, u_1,u_2 \}$.

In each case $A$ is an isolating set of both the $F_i$. Thus $x$ and all of $V(F_1) \cup V(F_2)$ become red.
Since it is incident with at least two edges, the weight of $x$ decreases by at least $\beta_2$.
Hence $\xi (A) \ge \expression_4 := |V(F_1) \cup V(F_2)|  \, \omega + \beta_2$. (We note that one can 
also consider the edges from the $F_i$ to $B\setminus \{x\}$, but this does not help the overall linear program.)
That is,
if $4\omega + \beta_2 \ge 1$, $7\omega + \beta_2 \ge 2$, and $10 \omega + \beta_2 \ge 3$, which system  $\cL$ ensures,
it holds that $\xi(A) \ge |A|$, so that $A$ is a $D$-desirable set in each case.~\smallqed

\begin{claim}
\label{c:final}
If $D$ is not an isolating set of $G$ then there is a $D$-desirable set.
\end{claim}
\proof 
By Claim~\ref{c:whiteComponent} 
we may assume that every component $F$ in $\residual [W]$ is either $K_2$ or $C_5$.
Let $B_F$ be the set of blue vertices adjacent to vertices in $F$.
By Claim~\ref{c:blueTwo} we may assume that if $z \in B_F$ then all its neighbors lie in $F$.

Assume first that $F = K_2$ with vertices $v_1$ and $v_2$, and let $A = \{v_1\}$. 
The two vertices in~$F$ as well as all the vertices in $B_F$ become red.
Since $F$ has order $2$, we have $B_F \subseteq B_1 \cup B_2$.
By the minimum degree, each vertex of $F$ is adjacent to at least $\delta-1$ blue vertices in $\residual$;
thus there are at least $ 2(\delta-1) $ edges joining $F$ to  $B_F$. 
The removal of an edge joining a vertex in $B_1$ to~$F$  decreases the weight by~$\beta_1$, 
while the collective removal of the two edges joining a vertex in $B_2$ to $F$ decreases the weight by~$\beta_2$. 
Hence 
$\xi (A) \ge \expression_{5a} := 2 \, \omega + 2 ( \delta-1)  \min ( \beta_1, \beta_2/2) $.
The system $\cL$ ensures $\expression_{5a} \ge 1$, and so
the set $A$ is a $D$-desirable set

Assume second that $F$ is the cycle $v_1v_2 \ldots v_{5}v_1$, and let $A = \{ v_2, v_5 \}$.
The five vertices in $F$ as well as all the vertices in $B_F$ become red.
By Claim~\ref{c:small} we may assume $B_4$ is empty.
By the minimum degree, each vertex of $F$ is adjacent to at least $\delta-2$ blue vertices in $\residual$;
thus there are at least $5( \delta-2)$ edges joining $F$ to  $B_F$. 
The removal of the $i$ edges joining a vertex in $B_i$ to $F$  decreases the weight by~$\beta_i$ in total.
Hence 
$\xi (A) \ge \expression_{5b}  := 5  \, \omega + 5 ( \delta - 2)  \min ( \beta_1, \beta_2/2, \beta_3/3) $.
The system $\cL$ ensures $\expression_{5b} \ge 2$, and so
the set $A$ is a $D$-desirable set.~\smallqed 

This completes the proof of Lemma~\ref{l:alluring}. We provide some consequences.

\begin{theorem}  \label{t:minFour}
If $G$ is a graph of order $n$ with minimum degree at least $4$, then $\iota(G) \le \frac{13}{41} n$.
\end{theorem}
\proof
The system $\cL$ is satisfied by setting
\[
   \omega= \frac{26}{82}   \qquad  \beta_4 = \frac{14}{82}  \qquad  \beta_3 = \frac{12}{82} \qquad \beta_2 = \frac{10}{82} \qquad \beta_1 = \frac{5}{82} .
\]
Thus by Lemma~\ref{l:alluring} the graph $G$ is alluring, and so $\iota(G) \le \omega n $ by Theorem~\ref{t:prototype}.~\QED

Unfortunately the technique is useless for $\delta=3$: the optimum for the linear program has $\omega > \frac{1}{3}$!
But for example for $\delta=5$, it gives an upper bound of $\frac{23}{78}n$.
Note that the $\delta=4$ and $\delta=5$ upper bounds are better than the original upper bounds provided in~\cite{CH}.


\section{Triangle-Free Graphs and Beyond}

Consider the system $\cL'$ where the constraints in $\cL$ for Claim~\ref{c:final} are replaced by:
\begin{quote}
\begin{tabular}{ll}
$2 \omega + 2 (\delta-1) \beta_1 \ge 1$ & (Claim 6) \\
$5  \, \omega + 5 ( \delta - 2)  \min ( \beta_1, \beta_2/2 ) \ge 2$ 
\end{tabular}
\end{quote}

\begin{lemma}  \label{l:alluringTfree}
Assume the weights satisfy the system of constraints $\cL'$ for some $\delta \ge 3$. Then every triangle-free
graph with minimum degree at least $\delta $ is alluring.
\end{lemma}

The proof of Lemma~\ref{l:alluringTfree} uses Claims~\ref{c:big} through \ref{c:blueTwo}.
At that point we may assume that: if $F$ is a component in $\residual [W]$ then it is either $K_2$ or $C_5$,
and any blue vertex adjacent to $F$ is adjacent to no other white component.

\begin{claim} 
\label{c:newFinal}
If the weights satisfy the system $\cL'$ and $D$ is not an isolating set of $G$, then there is a $D$-desirable set.
\end{claim}
\proof 
The proof is almost the same as that of Claim~\ref{c:final}. The difference is that,
by the lack of triangles, we know that when $F=K_2$ each blue vertex has at most one neighbor in $F$ and so $B_F \subseteq B_1$,
while when $F=C_5$ each blue vertex has at most two neighbors in $F$ and so $B_F \subseteq B_1 \cup B_2$.
Thus the constraints simplify to those given above.~\smallqed 

This completes the proof of Lemma~\ref{l:alluringTfree}. 

\begin{theorem}  \label{t:minFourTfree}
If $G$ is a triangle-free graph of order $n$ with minimum degree at least $4$, then $\iota(G) \le \frac{3}{10} n$.
\end{theorem}
\proof
The system $\cL'$ is satisfied by setting
\[
   \omega= \frac{3}{10}   \qquad  \beta_4 = \frac{3}{20}  \qquad  \beta_3 = \frac{1}{8} \qquad \beta_2 = \frac{1}{10} \qquad \beta_1 = \frac{1}{15} .
\]
The result follows as before.~\QED

For $\delta \ge 5$ the resultant bound is $\frac{9}{31} n$.
Unfortunately this system still does not provide a useful bound for $\delta=3$.
But if we add the requirement that the girth is at least five, then in Claim~\ref{c:final} one gets that for $F=C_5$ each blue vertex
has at most one neighbor in $F$, and so $B_F \subseteq B_1$. If we use use this fact to simplify the system of constraints,
then the optimum shows that: if $G$ is graph with $\delta \ge 3$ and girth at least $5$, then $\iota(G) \le \frac{11}{34}n$.
(Unfortunately, this is larger than the bound of~\cite{ZZ}, which holds for all cubic graphs.)

 \section{Lower Bounds}
 
 We noted in the introduction that Caro and Hansberg~\cite{CH} provided graphs of order $n$ with minimum degree $\delta$ where the isolation 
 number is $\frac{2}{\delta+4}n$. That family is not regular. For $\delta=3$ Ziemann and \.{Z}yli\'{n}ski~\cite{ZZ} found a family of regular graphs with 
 the same ratio $\frac{2}{7}$. 
 
 One can use the same approach as~\cite{ZZ} to generate some other regular examples. 
Namely,  find a regular graph $F$ of order $c$ that has a special edge $xy$ where 
each of $F$, $F-x$, $F-y$, and $F-\{x,y\}$ have isolation number (at least) $b$.
Then let $G_s$ be the graph formed by taking $s$ copies of the edge-deleted $F-xy$ 
and adding $s$ edges to make $G_s$ connected and $r$-regular. 
Every isolating set of the resultant $G_s$ must contain at least $b$ vertices from each copy of $F$: even if both $x$ and $y$ are dominated from outside,
one still needs $b$ vertices to isolate~$F$.

For $4$-regular, a suitable $F$ is for example the prism of $K_4$ 
with the special edge one joining the two copies of $K_4$. This example has $c=8$ and $b=2$;
so the associated $G_s$ has isolation number $\frac{1}{4}$ its order. 
For $4$-regular triangle-free graphs, there are several suitable graphs with $c=14$ and $b=3$.
One example is the metacirculant constructed by taking two $7$-cycles 
$u_1u_2 \ldots u_7u_1$ and $v_1v_2 \ldots v_7v_1$, and joining $u_i$ to both $v_{2i}$ and $v_{2i+3}$ (arithmetic modulo 7) for $i \in [7]$;
a special edge is one in the $u$-cycle. The associated $G_s$ has isolation number~$\frac{3}{14}n$ its order.
For $5$-regular, there are also many suitable graphs with $c=14$ and $b=3$. So the associated $G_s$ has isolation number~$\frac{3}{14}$ its order.

\section*{Acknowledgement}

We thank Geoff Boyer for some helpful discussions.
Research of Michael A. Henning was supported in part by the 
South African National Research Foundation, Grant Number 132588, and the University of Johannesburg.

\end{document}